\documentclass[11pt,a4paper]{amsart}
\usepackage[T1]{fontenc}
\usepackage{lmodern}
\usepackage{amsmath,amssymb,amsthm,amsaddr,mathtools,upgreek}
\usepackage[a4paper,margin=28mm]{geometry}
\usepackage{booktabs,microtype,xcolor,tikz,hyperref}
\usetikzlibrary{calc,angles,quotes,arrows.meta}
\definecolor{pedalblue}{RGB}{35,85,125}
\definecolor{pedalpurple}{RGB}{116,61,143}
\definecolor{pedalgray}{RGB}{106,118,125}
\hypersetup{colorlinks=true,linkcolor=pedalblue,citecolor=pedalpurple,urlcolor=pedalblue,
 pdftitle={Pedal Polygons Come Full Circle},pdfauthor={Michal Zwierzynski}}
\newtheorem{theorem}{Theorem}[section]
\newtheorem{proposition}[theorem]{Proposition}
\newtheorem{lemma}[theorem]{Lemma}
\newtheorem{corollary}[theorem]{Corollary}
\theoremstyle{remark}

\newcommand{\CC}{\mathbb C}
\newcommand{\RR}{\mathbb R}
\newcommand{\ii}{\mathrm i}
\newcommand{\ee}{\mathrm e}
\newcommand{\piconst}{\uppi}
\newcommand{\D}{\mathcal D}
\newcommand{\conv}{\operatorname{conv}}
\newcommand{\vertices}{\operatorname{vert}}
\newcommand{\dist}{\operatorname{dist}}
\newcommand{\doi}[1]{\href{https://doi.org/#1}{\nolinkurl{doi:#1}}}

\title[When Pedal Shapes Come Back: a Complex Story]{When Pedal Shapes Come Back: a Complex Story}
\author{Micha\l{} Zwierzy\'nski}
\email{Michal.Zwierzynski@pw.edu.pl}
\email{ORCID: 0000-0002-9627-1563}
\address{Warsaw University of~Technology\\
Faculty of~Mathematics and Information Science\\
ul. Koszykowa 75\\00-662 Warsaw, Poland}
\begin{document}
\raggedbottom
\begin{abstract}
Under the~usual nondegeneracy assumptions, repeatedly projecting a~fixed point onto the~side lines of~an~$n$-gon returns its shape after $n$ steps. We explain this classical theorem through a~short complex-number identity: the~multipliers of~one step move cyclically around the~polygon. The~same calculation allows a~different common angle of~projection at~each step. A~sharp angle-dependent contraction bound identifies the~incircle centre as~the~maximiser for~tangential polygons. In~space, even a~regular tetrahedron can fail to~recover its shape after four steps. We give an~exact example, then show how concentric inscribed and circumscribed spheres restore a~two-step return. Rectangular boxes provide a~further four-step family.
\end{abstract}
\maketitle

\section{A classical return}

\begin{figure}[tbp]
\centering
\begin{tikzpicture}[x=1.48cm,y=1.48cm,line join=round]
\begin{scope}[shift={(0.0,0.0)}]
\draw[pedalblue,line width=.8pt,fill=pedalblue!5] (0.000000000,0.874947662) -- (0.493079234,-0.848149878) -- (-0.963648329,0.267173911) -- (0.812508824,0.203931140) -- (-0.632963669,-0.707729487) -- cycle;
\fill[pedalblue] (0.000000000,0.874947662) circle (1.2pt);
\draw[pedalblue,line width=.45pt] (0.000000000,0.874947662) circle (2.8pt);
\fill[pedalblue] (0.493079234,-0.848149878) circle (1.2pt);
\fill[pedalblue] (-0.963648329,0.267173911) circle (1.2pt);
\fill[pedalblue] (0.812508824,0.203931140) circle (1.2pt);
\fill[pedalblue] (-0.632963669,-0.707729487) circle (1.2pt);
\fill[black] (0,0) circle (1pt);
\node[font=\small] at (0,-1.25) {$Z^{(0)}$};
\end{scope}
\begin{scope}[shift={(3.0,0.0)}]
\draw[pedalgray,line width=.8pt,fill=pedalgray!5] (-0.000000000,0.644164661) -- (0.596113928,-0.802899860) -- (-0.592254124,0.192528430) -- (0.670360664,0.195655505) -- (-0.532501031,-0.687313361) -- cycle;
\fill[pedalgray] (-0.000000000,0.644164661) circle (1.2pt);
\draw[pedalgray,line width=.45pt] (-0.000000000,0.644164661) circle (2.8pt);
\fill[pedalgray] (0.596113928,-0.802899860) circle (1.2pt);
\fill[pedalgray] (-0.592254124,0.192528430) circle (1.2pt);
\fill[pedalgray] (0.670360664,0.195655505) circle (1.2pt);
\fill[pedalgray] (-0.532501031,-0.687313361) circle (1.2pt);
\fill[black] (0,0) circle (1pt);
\node[font=\small] at (0,-1.25) {$Z^{(1)}$};
\end{scope}
\begin{scope}[shift={(6.0,0.0)}]
\draw[pedalgray,line width=.8pt,fill=pedalgray!5] (0.000000000,1.000000000) -- (0.440299945,-0.840074374) -- (-0.731801034,0.299345425) -- (0.945428219,0.233824831) -- (-0.679514676,-0.699074375) -- cycle;
\fill[pedalgray] (0.000000000,1.000000000) circle (1.2pt);
\draw[pedalgray,line width=.45pt] (0.000000000,1.000000000) circle (2.8pt);
\fill[pedalgray] (0.440299945,-0.840074374) circle (1.2pt);
\fill[pedalgray] (-0.731801034,0.299345425) circle (1.2pt);
\fill[pedalgray] (0.945428219,0.233824831) circle (1.2pt);
\fill[pedalgray] (-0.679514676,-0.699074375) circle (1.2pt);
\fill[black] (0,0) circle (1pt);
\node[font=\small] at (0,-1.25) {$Z^{(2)}$};
\end{scope}
\begin{scope}[shift={(0.0,-2.7)}]
\draw[pedalgray,line width=.8pt,fill=pedalgray!5] (0.000000000,0.626693843) -- (0.425780689,-0.672138732) -- (-0.701426504,0.197083305) -- (0.692183812,0.203769674) -- (-0.577216304,-0.816591292) -- cycle;
\fill[pedalgray] (0.000000000,0.626693843) circle (1.2pt);
\draw[pedalgray,line width=.45pt] (0.000000000,0.626693843) circle (2.8pt);
\fill[pedalgray] (0.425780689,-0.672138732) circle (1.2pt);
\fill[pedalgray] (-0.701426504,0.197083305) circle (1.2pt);
\fill[pedalgray] (0.692183812,0.203769674) circle (1.2pt);
\fill[pedalgray] (-0.577216304,-0.816591292) circle (1.2pt);
\fill[black] (0,0) circle (1pt);
\node[font=\small] at (0,-1.25) {$Z^{(3)}$};
\end{scope}
\begin{scope}[shift={(3.0,-2.7)}]
\draw[pedalgray,line width=.8pt,fill=pedalgray!5] (-0.000000000,0.710313782) -- (0.557002045,-0.819201570) -- (-0.694125061,0.223864259) -- (0.935796962,0.352539424) -- (-0.543686528,-0.649132718) -- cycle;
\fill[pedalgray] (-0.000000000,0.710313782) circle (1.2pt);
\draw[pedalgray,line width=.45pt] (-0.000000000,0.710313782) circle (2.8pt);
\fill[pedalgray] (0.557002045,-0.819201570) circle (1.2pt);
\fill[pedalgray] (-0.694125061,0.223864259) circle (1.2pt);
\fill[pedalgray] (0.935796962,0.352539424) circle (1.2pt);
\fill[pedalgray] (-0.543686528,-0.649132718) circle (1.2pt);
\fill[black] (0,0) circle (1pt);
\node[font=\small] at (0,-1.25) {$Z^{(4)}$};
\end{scope}
\begin{scope}[shift={(6.0,-2.7)}]
\draw[pedalpurple,line width=.8pt,fill=pedalpurple!5] (0.000000000,0.874947662) -- (0.493079234,-0.848149878) -- (-0.963648329,0.267173911) -- (0.812508824,0.203931140) -- (-0.632963669,-0.707729487) -- cycle;
\fill[pedalpurple] (0.000000000,0.874947662) circle (1.2pt);
\draw[pedalpurple,line width=.45pt] (0.000000000,0.874947662) circle (2.8pt);
\fill[pedalpurple] (0.493079234,-0.848149878) circle (1.2pt);
\fill[pedalpurple] (-0.963648329,0.267173911) circle (1.2pt);
\fill[pedalpurple] (0.812508824,0.203931140) circle (1.2pt);
\fill[pedalpurple] (-0.632963669,-0.707729487) circle (1.2pt);
\fill[black] (0,0) circle (1pt);
\node[font=\small] at (0,-1.25) {$Z^{(5)}$};
\end{scope}
\end{tikzpicture}
\caption{Five pedal steps for~a~self-intersecting pentagon shaped like an~irregular pentagram. Crossings of~sides are not vertices. Each panel is independently rotated about~$P$ (the~black dot) so that its marked first vertex lies on~the~positive vertical ray, and rescaled so that the~greatest distance from~$P$ is~$1$. The~first and last panels coincide}
\label{fig:pentagon}
\end{figure}
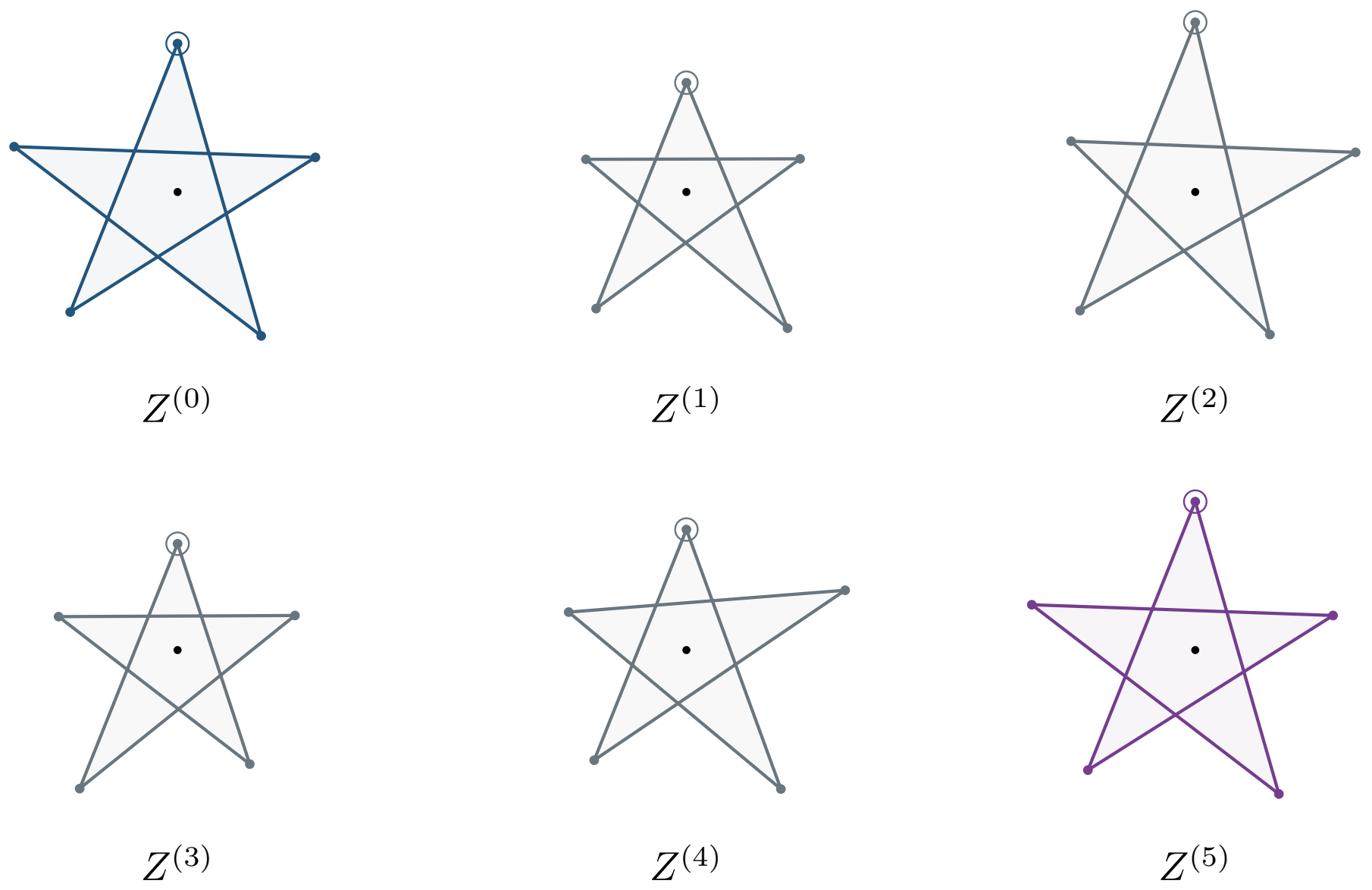

\noindent Choose a~point $P$ and a~triangle. Project $P$ perpendicularly onto the~three side lines, join the~feet, and repeat with the~same point. The~third triangle is similar to~the~first. The~intermediate triangles need not resemble either of~them. More generally, the~$n$th pedal $n$-gon is similar to~the~initial $n$-gon -- see Figure~\ref{fig:pentagon}.

The~triangle theorem appears as~Proposition~4 on~page~253 of~Casey's sixth edition~\cite{Casey}. Coxeter and Greitzer attribute its addition to~Neuberg~\cite[pp.~23--25]{CoxeterGreitzer}. The~polygon theorem also predates Stewart's 1940 paper. In~his historical note~\cite{Thebault1941}, Th\'ebault credits A.~Boutin with the~$n$-step similarity, J.~Dieudonn\'e with identifying the~fixed point as~its centre, and his own earlier work with the~distance-product formula for~its ratio. These attributions are taken from~that note.

Stewart treats the~more general iteration in~which the~common angle of~projection may change between steps~\cite{Stewart}. His proof follows two families of~directed angles: one stays fixed, while the~other shifts cyclically. De~Villiers later revisited the~construction through Miquel triangles~\cite{deVilliers}. We express the~same cyclic mechanism by~multiplication in~$\CC$, then compare the~planar return with its failure and some surviving special cases in~space.

Throughout the~planar discussion, a~polygon is a~closed polygonal chain given by~an~ordered cyclic list of~vertices. Neither convexity nor simplicity is required: self-intersections are allowed, as~Figure~\ref{fig:pentagon} illustrates. Its sides mean entire lines, so projections may lie on~extensions. We call an~$n$-gon \emph{admissible} if its first $n$ pedal steps are defined, with distinct consecutive vertices and with no side line used in~these steps passing through~$P$. These hypotheses avoid degenerate steps. Stewart's treatment also assigns side lines to~certain coincident vertices by~a~tangent construction~\cite[pp.~463--465]{Stewart}.

\section{Multipliers on the move}

\noindent Put $P=0$ in~the~complex plane. Multiplication by~a~nonzero complex number represents a~rotation and a~uniform scaling about~$P$. The~following identity is the~local reason for~the~return.

\begin{lemma}\label{lem:local}
Let $p$ and $q$ be distinct nonzero perpendicular projections of~$0$ onto two lines through $b\ne0$. The~perpendicular projection of~$0$ onto the~line $pq$ is
\begin{align*}
 r=\frac{pq}{b}.
\end{align*}
\end{lemma}

\begin{proof}
Perpendicularity means that $b-p$ is a~real multiple of~$\ii p$, and similarly for~$q$. Thus, for~some real numbers $s,t$,
\begin{align*}
 b=(1+\ii s)p=(1+\ii t)q.
\end{align*}
With $r=pq/b$, we obtain
\begin{align*}
 p=(1+\ii t)r,\qquad q=(1+\ii s)r.
\end{align*}
Hence $p-r=\ii t r$ and $q-p=\ii(s-t)r$. Since $p\ne q$, the~line $pq$ is precisely $r+\ii\RR r$. It passes through~$r$ and is perpendicular to~$0r$.
\end{proof}

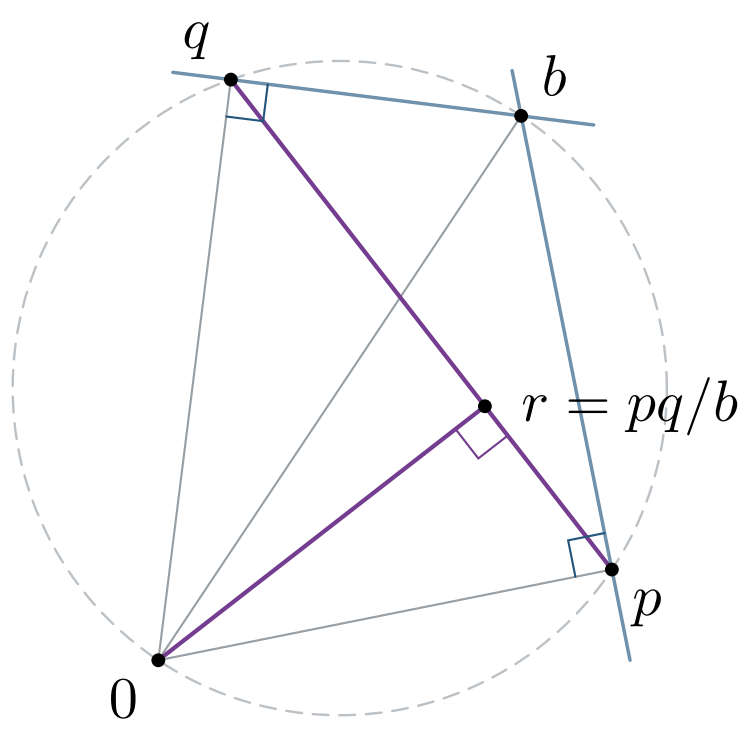
\begin{figure}[tbp]
\centering
\begin{tikzpicture}[x=1.2cm,y=1.2cm,line cap=round]
 \coordinate (O) at (0,0); \coordinate (b) at (2,3);
 \coordinate (p) at (2.5,.5); \coordinate (q) at (.4,3.2);
 \coordinate (r) at (1.8,1.4);
 \draw[pedalgray!45,densely dashed,line width=.45pt] (1,1.5) circle ({sqrt(3.25)});
 \draw[pedalblue!65,line width=.65pt] ($(p)!-.2!(b)$)--($(p)!1.1!(b)$);
 \draw[pedalblue!65,line width=.65pt] ($(q)!-.2!(b)$)--($(q)!1.25!(b)$);
 \draw[pedalgray!70] (O)--(p) (O)--(q) (O)--(b);
 \draw[pedalpurple,line width=.8pt] (p)--(q);
 \draw[pedalpurple,line width=.8pt] (O)--(r);
 \pic[draw=pedalblue,angle radius=7pt] {right angle=O--p--b};
 \pic[draw=pedalblue,angle radius=7pt] {right angle=O--q--b};
 \pic[draw=pedalpurple,angle radius=7pt] {right angle=O--r--p};
 \foreach \v in {O,b,p,q,r}{\fill (\v) circle (1.3pt);}
 \node[below left] at (O) {$0$}; \node[above right] at (b) {$b$};
 \node[below right] at (p) {$p$}; \node[above left] at (q) {$q$};
 \node[right=3pt] at (r) {$r=pq/b$};
\end{tikzpicture}
\caption{The~local identity. The~four points $0,p,b,q$ also lie on~a~circle of~diameter $0b$, but the~proof uses only multiplication by~$\ii$ to~encode perpendicularity}
\label{fig:local}
\end{figure}

Write $z_i^{(k)}$ for~the~$i$th vertex after $k$ pedal steps, with indices $i$ taken modulo~$n$, and let $\D$ denote one step. The~foot on~the~line $z_i^{(k)}z_{i+1}^{(k)}$ is labelled $z_i^{(k+1)}$. Applying Lemma~\ref{lem:local} with
\begin{align*}
 b=z_{i+1}^{(k)},\qquad p=z_i^{(k+1)},\qquad q=z_{i+1}^{(k+1)}
\end{align*}
gives the~shift identity
\begin{align}\label{eq:shift}
 \frac{z_i^{(k+2)}}{z_i^{(k+1)}}
 =\frac{z_{i+1}^{(k+1)}}{z_{i+1}^{(k)}}.
\end{align}
For~triangles, a~related real distance-ratio identity is used in~\cite[Figure~9 and equation~(1)]{HumenbergerDeVilliers}. Here the~complex quotient records the~rotation as~well as~the~scale.

\begin{theorem}\label{thm:return}
For~an~admissible $n$-gon $Z=(z_1,\ldots,z_n)$, put
\begin{align*}
 \lambda_i=\frac{z_i^{(1)}}{z_i},\qquad
 \Lambda=\prod_{i=1}^{n}\lambda_i.
\end{align*}
Then $\D^n(Z)=\Lambda Z$. In~particular, the~$n$th pedal polygon is directly similar to~$Z$, by~a~similarity fixing~$P$.
\end{theorem}

\begin{proof}
Equation~\eqref{eq:shift} moves the~multipliers one place at~each step. Inductively,
\begin{align*}
 z_i^{(m)}=\lambda_i\lambda_{i+1}\cdots\lambda_{i+m-1}z_i.
\end{align*}
For~$m=n$, every product contains each of~$\lambda_1,\ldots,\lambda_n$ exactly once. Complex multiplication is commutative, so the~product is independent of~$i$.
\end{proof}

The~scale can therefore be read from~the~first construction alone:
\begin{align}\label{eq:scale}
 |\Lambda|=\prod_{i=1}^{n}
 \frac{\dist(P,\operatorname{line}(A_iA_{i+1}))}{|PA_i|}.
\end{align}
This is the~classical distance-product formula recorded in~\cite[Theorem~3]{Thebault1941}, with the~reciprocal ratio convention there. It also appears for~triangles in~\cite[equation~(2)]{HumenbergerDeVilliers}. It lies strictly between $0$ and~$1$. Indeed, every factor is at~most~$1$. Equality of~all factors would make each side perpendicular to~its initial radius, giving $|z_{i+1}|^2=|z_i|^2+|z_{i+1}-z_i|^2$ around a~closed cycle, which is impossible. Moreover, $\D(cZ)=c\D(Z)$ for~$c\ne0$, so
\begin{align*}
 \D^{mn+j}(Z)=\Lambda^m\D^j(Z)
 \qquad(m\geq0,\ 0\leq j<n).
\end{align*}
Thus the~successive shapes repeat while the~vertices approach~$P$.

A~related inverse iteration is studied by~McDonald, Garcia and Reznik~\cite{McDonaldGarciaReznik}. Their circumcentre map replaces each pair $z_i,z_{i+1}$ by~the~circumcentre $c_i$ of~$0,z_i,z_{i+1}$. The~line $c_{i-1}c_i$ is the~perpendicular bisector of~$0z_i$, so its perpendicular foot is $z_i/2$. Thus the~pedal of~the~circumcentre polygon is $\tfrac12Z$, up to~a~cyclic shift of~labels, whenever these constructions are defined.

The~next estimate extends the~triangle bound in~\cite[Theorem~4]{HumenbergerDeVilliers}. It answers the~interior part of~the~question about tangential $n$-gons posed at~the~end of~\cite[Section~2.1]{HumenbergerDeVilliers}.

\begin{samepage}
\begin{corollary}\label{cor:bound}
Suppose $Z$ is admissible and strictly convex, with interior angles $\alpha_1,\ldots,\alpha_n$, and $P$ lies strictly inside it. Then
\begin{align}\label{eq:anglebound}
 |\Lambda|\leq\prod_{i=1}^{n}\sin\frac{\alpha_i}{2}
 \leq\cos^n\frac{\piconst}{n}.
\end{align}
Equality in~the~first inequality holds exactly at~an~incircle centre. Equality with the~rightmost bound holds exactly for~a~regular polygon projected from~its centre. For~a~tangential polygon, the~incircle centre is admissible and uniquely maximises the~scale and the~$n$th pedal's area over admissible interior points.
\end{corollary}
\end{samepage}

\begin{proof}
Put $r_i=|PA_i|$ and $d_i=\dist(P,\operatorname{line}(A_iA_{i+1}))$. The~segment $A_iP$ splits the~interior angle into $\beta_i,\gamma_i>0$, with $\beta_i+\gamma_i=\alpha_i$. Thus
\begin{align*}
 |\Lambda|^2
 =\prod_{i=1}^{n}\frac{d_{i-1}d_i}{r_i^2}
 =\prod_{i=1}^{n}\sin\beta_i\sin\gamma_i
 \leq\prod_{i=1}^{n}\sin^2\frac{\alpha_i}{2}.
\end{align*}
Here $\sin\beta_i\sin\gamma_i\leq\sin^2(\alpha_i/2)$, with equality exactly when $\beta_i=\gamma_i$. Equality in~the~product means that $P$ lies on~every internal angle bisector, equivalently $d_1=\cdots=d_n$: the~incircle condition. Concavity of~$\log\sin x$ and $\sum_i\alpha_i=(n-2)\piconst$ give the~second inequality, with equality only for~equal angles. An~equiangular tangential polygon is regular.

For~completeness, an~incircle centre is admissible. Put it at~$0$, orient the~polygon counterclockwise, and set $\varepsilon_i=\piconst-\alpha_i\in(0,\piconst)$. The~first-step multipliers are
\begin{align*}
 \lambda_i=\cos\frac{\varepsilon_i}{2}\,
            \ee^{\ii\varepsilon_i/2}.
\end{align*}
The~candidate vertices $z_i\prod_{j=0}^{m-1}\lambda_{i+j}$ have successive polar increments $(\varepsilon_{i+1}+\varepsilon_{i+m})/2\in(0,\piconst)$. They are nonzero, consecutive vertices are distinct, and no side line passes through~$0$. Lemma~\ref{lem:local} therefore identifies these candidates inductively with the~actual iterates. The~first bound is attained at~this point and nowhere else -- the~area ratio is $|\Lambda|^2$.
\end{proof}

Stewart's theorem also allows the~angle of~projection to~change at~every step. To~state this explicitly, fix $\theta\in(0,\piconst)$ and replace each side line $\ell$ by~the~point $q\in\ell$ for~which the~directed angle from~the~line $0q$ to~$\ell$, measured counterclockwise modulo~$\piconst$, is~$\theta$. Denote the~resulting polygon by~$M_\theta(Z)$. Thus $M_{\piconst/2}=\D$.

\begin{corollary}\label{cor:angles}
Let $Z$ be an~admissible $n$-gon and choose arbitrary angles $\theta_1,\ldots,\theta_n\in(0,\piconst)$. Keep the~point $P$ fixed and use the~common angle $\theta_k$ for~all side lines in~step~$k$. Then
\begin{align}\label{eq:oblique}
 M_{\theta_n}\circ\cdots\circ M_{\theta_1}(Z)
 =\Lambda\prod_{k=1}^{n}(1+\ii\cot\theta_k)\,Z.
\end{align}
In~particular, the~$n$th polygon is directly similar to~$Z$, by~a~similarity fixing~$P$, even when the~angles in~different steps are different.
\end{corollary}

\begin{proof}
A~side with perpendicular foot $p\ne0$ is the~line $p+\ii\RR p$. The~point $q=(1+\ii\cot\theta)p$ lies on~this line and has the~required angle, since
\begin{align*}
 1+\ii\cot\theta
 =\frac{\ee^{\ii(\piconst/2-\theta)}}{\sin\theta}.
\end{align*}
Consequently $M_\theta(Z)=(1+\ii\cot\theta)\D(Z)$. The~one-step similarity as~the~angle varies was already observed by~Th\'ebault~\cite[pp.~433--434]{Thebault1914}. As~$\D(cZ)=c\D(Z)$ for~$c\ne0$, the~scalar factors multiply, and every intermediate polygon is a~nonzero multiple of~its perpendicular pedal counterpart. Admissibility is therefore preserved, and Theorem~\ref{thm:return} gives~\eqref{eq:oblique}.
\end{proof}

\begin{figure}[tbp]
\centering
\begingroup
\definecolor{pedalblue}{RGB}{35,85,125}
\definecolor{pedalpurple}{RGB}{116,61,143}
\definecolor{pedalgray}{RGB}{112,119,128}
\begin{tikzpicture}[x=1cm,y=1cm,font=\small,line cap=round,line join=round]
\begin{scope}[xshift=0cm,yshift=0cm]
\node at (0,2.55) {\textcolor{pedalblue}{$Z^{(0)}$} $\longmapsto$ \textcolor{pedalpurple}{$Z^{(1)}$} $\qquad\theta_{1}=60^\circ$};
\fill[pedalblue!3] (-2.200000000,-1.350000000)--(2.200000000,-1.350000000)--(0.550000000,1.950000000)--cycle;
\fill[pedalpurple!5] (0.635085296,-1.350000000)--(0.938948822,1.172102355)--(-1.121770252,-0.056124302)--cycle;
\draw[pedalblue,line width=.8pt] (-2.200000000,-1.350000000)--(2.200000000,-1.350000000)--(0.550000000,1.950000000)--cycle;
\draw[pedalgray,densely dashed,line width=.4pt] (0.000000000,-0.250000000)--(0.635085296,-1.350000000);
\draw[pedalgray,densely dashed,line width=.4pt] (0.000000000,-0.250000000)--(0.938948822,1.172102355);
\draw[pedalgray,densely dashed,line width=.4pt] (0.000000000,-0.250000000)--(-1.121770252,-0.056124302);
\draw[pedalpurple,line width=.85pt] (0.635085296,-1.350000000)--(0.938948822,1.172102355)--(-1.121770252,-0.056124302)--cycle;
\draw[pedalgray,line width=.4pt,->] (0.517585296,-1.146484030) arc[start angle=120.0000000,delta angle=60,radius=0.235cm];
\node[font=\tiny,inner sep=.5pt,fill=white,text=pedalgray] at (0.262694372,-1.135000000) {$60^\circ$};
\draw[pedalgray,line width=.4pt,->] (0.809466203,0.975992052) arc[start angle=-123.4349488,delta angle=60,radius=0.235cm];
\node[font=\tiny,inner sep=.5pt,fill=white,text=pedalgray] at (0.913185260,0.742874865) {$60^\circ$};
\draw[pedalgray,line width=.4pt,->] (-0.890203286,-0.096146051) arc[start angle=-9.8055711,delta angle=60,radius=0.235cm];
\node[font=\tiny,inner sep=.5pt,fill=white,text=pedalgray] at (-0.718203817,0.092314684) {$60^\circ$};
\fill[pedalblue] (-2.200000000,-1.350000000) circle[radius=1.05pt];
\node[font=\scriptsize,anchor=north east,inner sep=2pt,text=pedalblue] at (-2.200000000,-1.350000000) {$z_{1}^{(0)}$};
\fill[pedalblue] (2.200000000,-1.350000000) circle[radius=1.05pt];
\node[font=\scriptsize,anchor=north west,inner sep=2pt,text=pedalblue] at (2.200000000,-1.350000000) {$z_{2}^{(0)}$};
\fill[pedalblue] (0.550000000,1.950000000) circle[radius=1.05pt];
\node[font=\scriptsize,anchor=south,inner sep=2pt,text=pedalblue] at (0.550000000,1.950000000) {$z_{3}^{(0)}$};
\fill[pedalpurple] (0.635085296,-1.350000000) circle[radius=1.3pt];
\fill[pedalpurple] (0.938948822,1.172102355) circle[radius=1.3pt];
\fill[pedalpurple] (-1.121770252,-0.056124302) circle[radius=1.3pt];
\node[font=\scriptsize,anchor=north,inner sep=2pt,text=pedalpurple] at (0.635085296,-1.350000000) {$z_1^{(1)}$};
\fill[pedalgray] (0.000000000,-0.250000000) circle[radius=1pt];
\node[font=\scriptsize,anchor=north east,inner sep=2pt,fill=white] at (0.000000000,-0.250000000) {$P$};
\end{scope}
\begin{scope}[xshift=7.0cm,yshift=0cm]
\node at (0,2.55) {\textcolor{pedalblue}{$Z^{(1)}$} $\longmapsto$ \textcolor{pedalpurple}{$Z^{(2)}$} $\qquad\theta_{2}=90^\circ$};
\fill[pedalblue!3] (-2.200000000,-1.350000000)--(2.200000000,-1.350000000)--(-0.339016393,1.939180328)--cycle;
\fill[pedalpurple!5] (-0.440000000,-1.350000000)--(0.575776278,0.754108004)--(-1.207475410,0.404229508)--cycle;
\draw[pedalblue,line width=.8pt] (-2.200000000,-1.350000000)--(2.200000000,-1.350000000)--(-0.339016393,1.939180328)--cycle;
\draw[pedalgray,densely dashed,line width=.4pt] (-0.440000000,-0.030000000)--(-0.440000000,-1.350000000);
\draw[pedalgray,densely dashed,line width=.4pt] (-0.440000000,-0.030000000)--(0.575776278,0.754108004);
\draw[pedalgray,densely dashed,line width=.4pt] (-0.440000000,-0.030000000)--(-1.207475410,0.404229508);
\draw[pedalpurple,line width=.85pt] (-0.440000000,-1.350000000)--(0.575776278,0.754108004)--(-1.207475410,0.404229508)--cycle;
\draw[pedalgray,line width=.4pt] (-0.440000000,-1.115000000)--(-0.675000000,-1.115000000)--(-0.675000000,-1.350000000);
\node[font=\tiny,inner sep=.5pt,fill=white,text=pedalgray] at (-0.744055916,-1.045944084) {$90^\circ$};
\draw[pedalgray,line width=.4pt] (0.389752552,0.610510742)--(0.533349814,0.424487017)--(0.719373539,0.568084278);
\node[font=\tiny,inner sep=.5pt,fill=white,text=pedalgray] at (0.520882587,0.327626254) {$90^\circ$};
\draw[pedalgray,line width=.4pt] (-1.002943285,0.288507385)--(-0.887221161,0.493039510)--(-1.091753286,0.608761633);
\node[font=\tiny,inner sep=.5pt,fill=white,text=pedalgray] at (-0.793112862,0.519136769) {$90^\circ$};
\fill[pedalblue] (-2.200000000,-1.350000000) circle[radius=1.05pt];
\node[font=\scriptsize,anchor=north east,inner sep=2pt,text=pedalblue] at (-2.200000000,-1.350000000) {$z_{1}^{(1)}$};
\fill[pedalblue] (2.200000000,-1.350000000) circle[radius=1.05pt];
\node[font=\scriptsize,anchor=north west,inner sep=2pt,text=pedalblue] at (2.200000000,-1.350000000) {$z_{2}^{(1)}$};
\fill[pedalblue] (-0.339016393,1.939180328) circle[radius=1.05pt];
\node[font=\scriptsize,anchor=south,inner sep=2pt,text=pedalblue] at (-0.339016393,1.939180328) {$z_{3}^{(1)}$};
\fill[pedalpurple] (-0.440000000,-1.350000000) circle[radius=1.3pt];
\fill[pedalpurple] (0.575776278,0.754108004) circle[radius=1.3pt];
\fill[pedalpurple] (-1.207475410,0.404229508) circle[radius=1.3pt];
\node[font=\scriptsize,anchor=north,inner sep=2pt,text=pedalpurple] at (-0.440000000,-1.350000000) {$z_1^{(2)}$};
\fill[pedalgray] (-0.440000000,-0.030000000) circle[radius=1pt];
\node[font=\scriptsize,anchor=north east,inner sep=2pt,fill=white] at (-0.440000000,-0.030000000) {$P$};
\end{scope}
\begin{scope}[xshift=0cm,yshift=-5.1cm]
\node at (0,2.55) {\textcolor{pedalblue}{$Z^{(2)}$} $\longmapsto$ \textcolor{pedalpurple}{$Z^{(3)}$} $\qquad\theta_{3}=120^\circ$};
\fill[pedalblue!3] (-2.200000000,-1.350000000)--(2.200000000,-1.350000000)--(0.146666667,1.387777778)--cycle;
\fill[pedalpurple!5] (-0.585346958,-1.350000000)--(1.277523964,-0.120031952)--(-0.343528668,0.815883221)--cycle;
\draw[pedalblue,line width=.8pt] (-2.200000000,-1.350000000)--(2.200000000,-1.350000000)--(0.146666667,1.387777778)--cycle;
\draw[pedalgray,densely dashed,line width=.4pt] (0.038596491,-0.269298246)--(-0.585346958,-1.350000000);
\draw[pedalgray,densely dashed,line width=.4pt] (0.038596491,-0.269298246)--(1.277523964,-0.120031952);
\draw[pedalgray,densely dashed,line width=.4pt] (0.038596491,-0.269298246)--(-0.343528668,0.815883221);
\draw[pedalpurple,line width=.85pt] (-0.585346958,-1.350000000)--(1.277523964,-0.120031952)--(-0.343528668,0.815883221)--cycle;
\draw[pedalgray,line width=.4pt,->] (-0.467846958,-1.146484030) arc[start angle=60.0000000,delta angle=120,radius=0.235cm];
\node[font=\tiny,inner sep=.5pt,fill=white,text=pedalgray] at (-0.800346958,-0.977609076) {$120^\circ$};
\draw[pedalgray,line width=.4pt,->] (1.044211188,-0.148141534) arc[start angle=-173.1301024,delta angle=120,radius=0.235cm];
\node[font=\tiny,inner sep=.5pt,fill=white,text=pedalgray] at (1.108611225,-0.515466506) {$120^\circ$};
\draw[pedalgray,line width=.4pt,->] (-0.265475810,0.594224133) arc[start angle=-70.6012946,delta angle=120,radius=0.235cm];
\node[font=\tiny,inner sep=.5pt,fill=white,text=pedalgray] at (0.079131745,0.736774590) {$120^\circ$};
\fill[pedalblue] (-2.200000000,-1.350000000) circle[radius=1.05pt];
\node[font=\scriptsize,anchor=north east,inner sep=2pt,text=pedalblue] at (-2.200000000,-1.350000000) {$z_{1}^{(2)}$};
\fill[pedalblue] (2.200000000,-1.350000000) circle[radius=1.05pt];
\node[font=\scriptsize,anchor=north west,inner sep=2pt,text=pedalblue] at (2.200000000,-1.350000000) {$z_{2}^{(2)}$};
\fill[pedalblue] (0.146666667,1.387777778) circle[radius=1.05pt];
\node[font=\scriptsize,anchor=south,inner sep=2pt,text=pedalblue] at (0.146666667,1.387777778) {$z_{3}^{(2)}$};
\fill[pedalpurple] (-0.585346958,-1.350000000) circle[radius=1.3pt];
\fill[pedalpurple] (1.277523964,-0.120031952) circle[radius=1.3pt];
\fill[pedalpurple] (-0.343528668,0.815883221) circle[radius=1.3pt];
\node[font=\scriptsize,anchor=north,inner sep=2pt,text=pedalpurple] at (-0.585346958,-1.350000000) {$z_1^{(3)}$};
\fill[pedalgray] (0.038596491,-0.269298246) circle[radius=1pt];
\node[font=\scriptsize,anchor=north east,inner sep=2pt,fill=white] at (0.038596491,-0.269298246) {$P$};
\end{scope}
\begin{scope}[xshift=7cm,yshift=-5.1cm]
\node at (0,2.55) {\textcolor{pedalblue}{$Z^{(0)}$} and \textcolor{pedalpurple}{$\mu^{-1}Z^{(3)}$}};
\fill[pedalblue!3] (-2.200000000,-1.350000000)--(2.200000000,-1.350000000)--(0.550000000,1.950000000)--cycle;
\draw[pedalblue,line width=1.35pt] (-2.200000000,-1.350000000)--(2.200000000,-1.350000000)--(0.550000000,1.950000000)--cycle;
\draw[pedalpurple,dash pattern=on 3.5pt off 3.5pt,line width=.75pt] (-2.200000000,-1.350000000)--(2.200000000,-1.350000000)--(0.550000000,1.950000000)--cycle;
\draw[pedalpurple,fill=white,line width=.65pt] (-2.200000000,-1.350000000) circle[radius=2.6pt];
\fill[pedalblue] (-2.200000000,-1.350000000) circle[radius=1pt];
\node[font=\scriptsize,anchor=north east,inner sep=3pt,text=pedalblue] at (-2.200000000,-1.350000000) {$z_{1}^{(0)}$};
\draw[pedalpurple,fill=white,line width=.65pt] (2.200000000,-1.350000000) circle[radius=2.6pt];
\fill[pedalblue] (2.200000000,-1.350000000) circle[radius=1pt];
\node[font=\scriptsize,anchor=north west,inner sep=3pt,text=pedalblue] at (2.200000000,-1.350000000) {$z_{2}^{(0)}$};
\draw[pedalpurple,fill=white,line width=.65pt] (0.550000000,1.950000000) circle[radius=2.6pt];
\fill[pedalblue] (0.550000000,1.950000000) circle[radius=1pt];
\node[font=\scriptsize,anchor=south,inner sep=3pt,text=pedalblue] at (0.550000000,1.950000000) {$z_{3}^{(0)}$};
\fill[pedalgray] (0.000000000,-0.250000000) circle[radius=1pt];
\node[font=\scriptsize,anchor=north east,inner sep=2pt] at (0.000000000,-0.250000000) {$P$};
\node[font=\scriptsize] at (0,-2.0) {$\displaystyle\mu=-\frac{4032+56\,\mathrm{i}}{25925}$};
\end{scope}
\end{tikzpicture}
\endgroup

\caption{Three oblique steps at~$60^\circ$, $90^\circ$, and $120^\circ$. Blue: the~input triangle; purple: the~new triangle; dashed: the~projection segments. Each pair and~$P$ are jointly rotated and rescaled for~legibility. The~marked angles run counterclockwise from~the~ray $qP$ to~the~backward side ray. All nine constructed points lie inside their sides. The~last panel overlays $Z^{(0)}$ with $\mu^{-1}Z^{(3)}$
}
\label{fig:angles}
\end{figure}
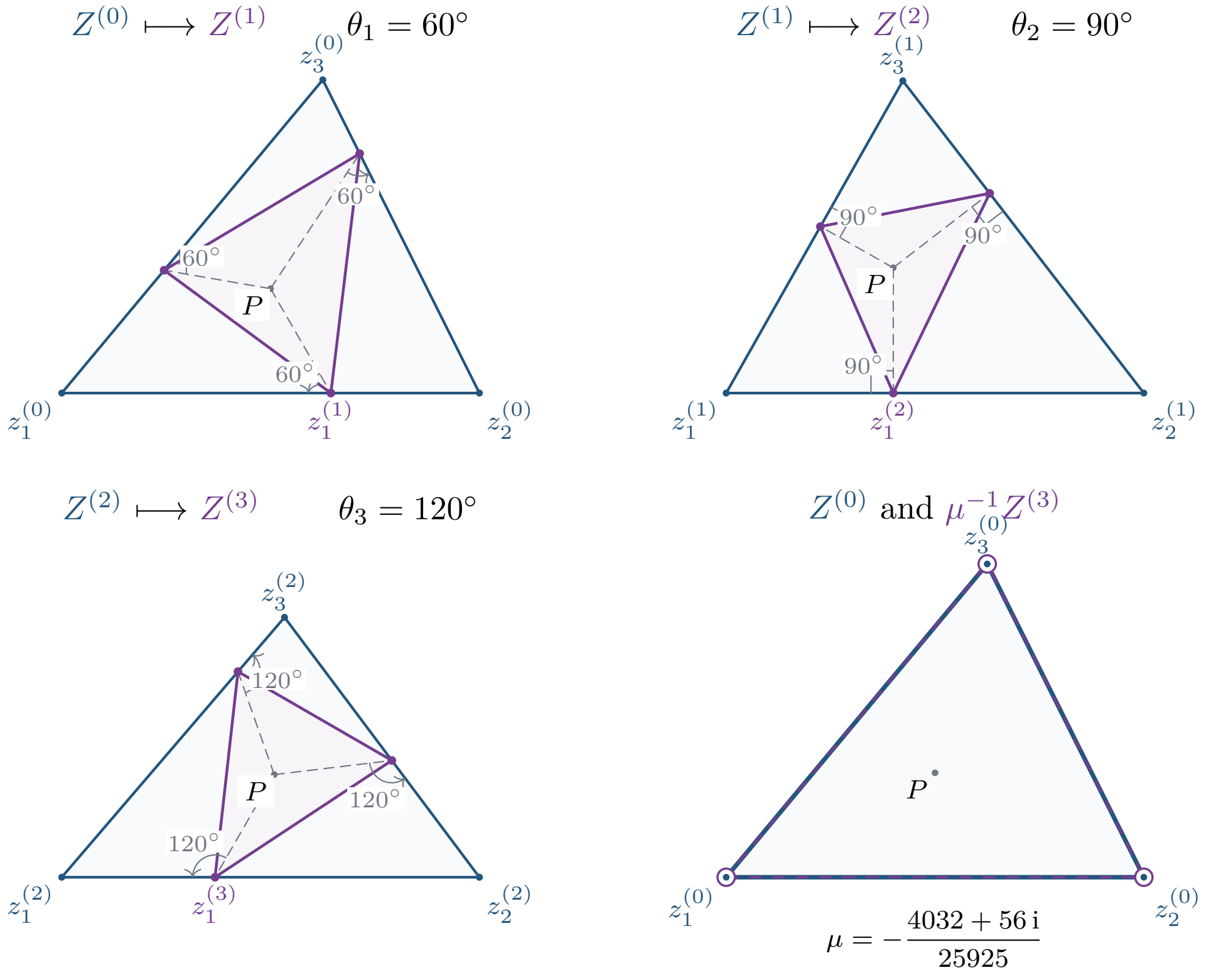

The~angles are common to~the~sides within one step. They are independent between steps. Formula~\eqref{eq:oblique} also shows that reordering the~chosen angles leaves the~final polygon unchanged. Its scale factor is
\begin{align*}
 \frac{|\Lambda|}{\prod_{k=1}^{n}\sin\theta_k},
\end{align*}
which need not be less than~$1$. For~fixed projection angles, Corollary~\ref{cor:bound} still makes the~incircle centre the~unique interior maximiser for~a~tangential polygon. For~the~angles $60^\circ,90^\circ,120^\circ$, the~product in~\eqref{eq:oblique} equals $4/3$. Figure~\ref{fig:angles} shows these steps for~$Z^{(0)}=(-2-\ii,2-\ii,\tfrac12+2\ii)$ with $P=0$.

\section{A regular tetrahedron that does not return}

\noindent For~a~convex polytope $K\subset\RR^d$ containing $P$ in~its interior, define $\D_P(K)$ to~be the~convex hull of~the~perpendicular projections of~$P$ onto the~hyperplanes containing its facets. For~a~tetrahedron these four feet form its pedal tetrahedron. We again put $P=0$ and write~$\D$.

It is tempting to~expect a~four-step return in~dimension~$3$. Such an~assertion appears in~\cite{Kostanyan}. The~following example disproves it even when the~initial tetrahedron is regular and every projection falls inside the~corresponding triangular face.

Define
\begin{align*}
 V(a,b,c)=\conv\{(a,b,b),(b,a,b),(b,b,a),(c,c,c)\}.
\end{align*}
The~first three vertices form an~equilateral base. Its three edges have squared length $B^2=2(a-b)^2$, while the~three lateral edges have squared length
\begin{align*}
 L^2=(a-c)^2+2(b-c)^2.
\end{align*}
Take $T_0=V(3,-2,3)$. All its edges have squared length~$50$, and
\begin{align*}
 0=\tfrac3{10}(3,-2,-2)+\tfrac3{10}(-2,3,-2)
   +\tfrac3{10}(-2,-2,3)+\tfrac1{10}(3,3,3),
\end{align*}
so $P$ is strictly interior.

The~calculation of~four iterates stays within this family. The~face opposite $(a,b,b)$ has equation
\begin{align*}
 (a+b-2c)x+(c-b)y+(c-b)z=c(a-b).
\end{align*}
Projecting~$0$ onto this plane, and onto the~base plane $x+y+z=a+2b$, gives
\begin{align*}
 \D(V(a,b,c))&=V\left(
 \frac{c(a-b)(a+b-2c)}{\Delta},
 \frac{c(a-b)(c-b)}{\Delta},
 \frac{a+2b}{3}\right),\\[-.5ex]
 \Delta&=(a+b-2c)^2+2(c-b)^2.
\end{align*}
Table~\ref{tab:tetra} records the~exact results using integer representatives. Negative scale factors mean central reflections as~well as~rescalings.

\begin{table}[!htbp]
\centering
\caption{Exact iterates: $\D^j(T_0)=\sigma_jV(a_j,b_j,c_j)$}
\label{tab:tetra}
\begin{tabular}{@{}crrrcc@{}}
\toprule
$j$&$a_j$&$b_j$&$c_j$&$\sigma_j$&$L_j^2/B_j^2$\\
\midrule
0&3&-2&3&$1$&$1$\\
1&3&-3&1&$-1/3$&$1/2$\\
2&1&-2&3&$1/9$&$3$\\
3&7&-5&11&$-1/99$&$11/6$\\
4&55&-44&19&$1/1881$&$57/121$\\
\bottomrule
\end{tabular}
\end{table}

\begin{figure}[tbp]
\centering
\begingroup
\definecolor{pedalblue}{RGB}{35,85,125}
\definecolor{pedalpurple}{RGB}{116,61,143}
\definecolor{pedalgray}{RGB}{112,119,128}
\begin{tikzpicture}[font=\small,
    x={(2.3cm,0cm)},y={(0cm,0.78664633cm)},
    z={(0cm,2.16129303cm)},
    line cap=round,line join=round]
  \foreach \shift/\hh/\pp/\col/\tag/\ratio in {
    0/1.4142135624/0.1414213562/pedalblue/0/{\sqrt{2/3}},
    6.4/0.6428243465/0.2357022604/pedalpurple/4/{5\sqrt{6}/33}}
  {
    \begin{scope}[xshift=\shift cm]
      \coordinate (A) at (0.2588190451,-0.9659258263,0);
      \coordinate (B) at (0.7071067812,0.7071067812,0);
      \coordinate (C) at (-0.9659258263,0.2588190451,0);
      \coordinate (S) at (0,0,\hh);
      \coordinate (G) at (0,0,0);
      \coordinate (P) at (0,0,\pp);
      \fill[\col!7] (S)--(A)--(B)--cycle;
      \fill[\col!16] (S)--(C)--(A)--cycle;
      \draw[pedalgray,densely dashed,line width=.45pt]
        (B)--(C);
      \draw[\col,line width=.85pt]
        (S)--(B)--(A)--(C)--(S) (S)--(A);
      \draw[pedalgray,densely dashed,line width=.4pt]
        (G)--(S);
      \fill[\col] (P) circle[radius=1.35pt];
      \node[anchor=east,inner sep=1pt,xshift=-2pt,
        fill=\col!16,text=\col]
        at (P) {$P$};
      \node[fill=white,inner sep=1.3pt,text=\col]
        at (barycentric cs:A=.5,B=.5) {$b$};
      \node[fill=white,inner sep=1.3pt,text=\col]
        at (barycentric cs:S=.53,C=.47) {$\ell$};
      \node[text=\col] at (0,0,1.67) {$T_{\tag}$};
      \node at (0,0,-.67) {$h/b=\ratio$};
    \end{scope}
  }
\end{tikzpicture}
\endgroup

\caption{The~initial regular tetrahedron and its fourth pedal, independently rescaled to~the~same base edge length~$b$. Here $h$ is the~height above the~base and $P$ is the~projection point. The~common viewing direction makes the~change of~shape visible}
\label{fig:tetra}
\end{figure}
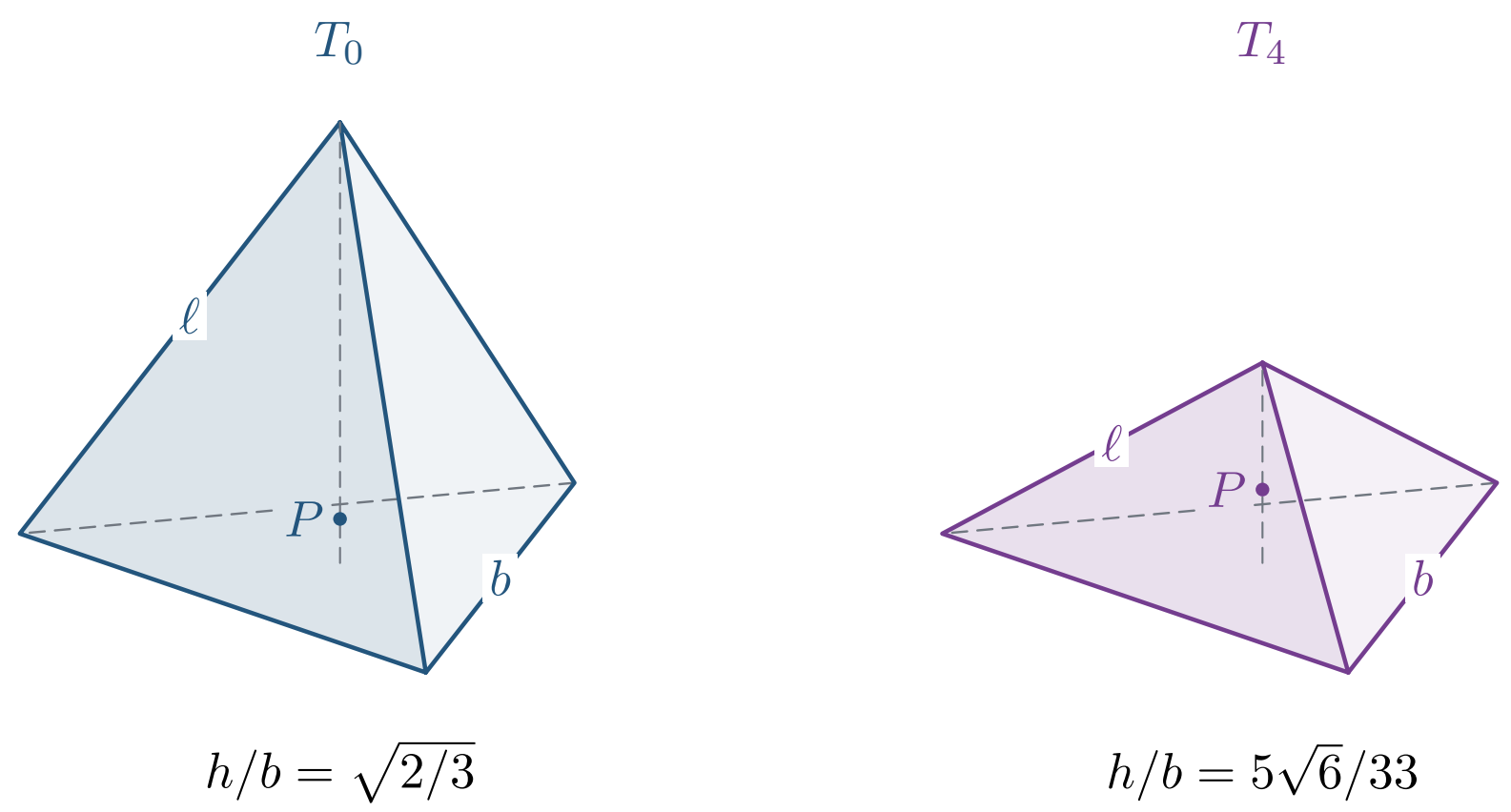

In~particular,
\begin{align*}
 \D^4(T_0)=\frac1{1881}V(55,-44,19),
 \qquad \frac{L_4^2}{B_4^2}=\frac{57}{121}\ne1.
\end{align*}
The~fourth tetrahedron is not regular, hence cannot be similar to~$T_0$, under any relabelling of~its vertices. Figure~\ref{fig:tetra} compares the~two shapes.

The~example does not exploit projections onto face extensions. For~the~three equivalent lateral faces, the~barycentric coordinates of~the~feet at~the~four steps are, up to~ordering,
\begin{align*}
 (\tfrac25,\tfrac25,\tfrac15),\quad
 (\tfrac16,\tfrac16,\tfrac23),\quad
 (\tfrac4{11},\tfrac4{11},\tfrac3{11}),\quad
 (\tfrac{33}{76},\tfrac{33}{76},\tfrac5{38}).
\end{align*}
The~base foot is always its centroid. All coordinates are positive. The~origin also stays interior: for~every row of~the~table, $a_j+2b_j<0<c_j$, placing it strictly between the~base centre and apex.

The~planar proof depends on~a~single complex multiplier and its cyclic transport. Spatial rotations do not in~general commute, and the~right-angle configuration behind Figure~\ref{fig:local} need not be coplanar in~space. These observations explain why that proof does not transfer. The~explicit calculation establishes the~failure of~the~four-step assertion.

\section{Two spheres, and a family of boxes}

\noindent There is nevertheless a~simple positive spatial statement. The~hypothesis involves \emph{both} spheres: every facet is tangent to~one sphere, and every vertex lies on~another, with the~same centre~$P$. The~polytope need not be regular.

\begin{proposition}\label{prop:spheres}
Let $K\subset\RR^d$ be a~full-dimensional convex polytope circumscribed about the~sphere of~radius~$r$ centred at~$0$. Then
\begin{align}\label{eq:inversion}
 \D^2(K)=\conv\left\{\frac{r^2v}{|v|^2}:v\in\vertices K\right\}.
\end{align}
In~particular, if all vertices have distance~$R$ from~$0$, then
\begin{align}\label{eq:spheres}
 \D^2(K)=\frac{r^2}{R^2}K.
\end{align}
Conversely, within this circumscribed class, a~positive homothety $\D^2(K)=\lambda K$ centred at~$0$ is possible only when all vertices have the~same distance from~$0$.
\end{proposition}

\begin{proof}
Write the~facet inequalities as~$u_F\cdot x\leq r$, where $|u_F|=1$. The~feet are $ru_F$. Recall that the~polar polytope is
\begin{align*}
 K^\circ=\{y:x\cdot y\leq1\text{ for~every }x\in K\}
          =\conv\{u_F/r:F\text{ a~facet}\}.
\end{align*}
Thus $\D(K)=r^2K^\circ$. By~the~standard vertex--facet correspondence for~polar polytopes~\cite[Chapter~2]{Ziegler}, its facets correspond to~the~vertices $v$ of~$K$, and their supporting hyperplanes are $v\cdot y=r^2$. The~projection of~$0$ onto such a~hyperplane is $r^2v/|v|^2$, proving~\eqref{eq:inversion} and~\eqref{eq:spheres}.

For~the~converse, if $\D^2(K)=\lambda K$ with $\lambda>0$, each vertex $\lambda v$ must equal $r^2w/|w|^2$ for~some vertex~$w$ of~$K$. The~vertices $v,w$ lie on~the~same positive ray. A~convex body containing~$0$ in~its interior has only one boundary point on~that ray, so $v=w$ and $|v|^2=r^2/\lambda$.
\end{proof}

For~example, the~tetrahedron with vertices
\begin{align*}
 (a,b,c),\quad(a,-b,-c),\quad(-a,b,-c),\quad(-a,-b,c),
 \qquad a,b,c>0,
\end{align*}
has concentric spheres with $R^2=a^2+b^2+c^2$ and $r^{-2}=a^{-2}+b^{-2}+c^{-2}$. Its second pedal is therefore homothetic to~it, even when it is not regular. The~usual regular convex polyhedra are other immediate examples.

A~different return survives for~rectangular boxes, which generally do not have a~concentric insphere tangent to~every facet.

\begin{samepage}
\begin{proposition}\label{prop:boxes}
Let $B=\prod_{i=1}^{d}[-a_i,a_i]$, where $a_i>0$, and put
\begin{align*}
 S=\sum_{i=1}^{d}a_i^{-2},\qquad T=\sum_{i=1}^{d}a_i^2.
\end{align*}
With projections taken from~the~centre~$0$,
\begin{align*}
 \D^4(B)=\frac1{ST}B.
\end{align*}
\end{proposition}
\end{samepage}

\begin{proof}
The~first pedal is the~cross-polytope $\conv\{\pm a_i e_i\}$, whose facet hyperplanes are
\begin{align*}
 \sum_{i=1}^{d}\frac{\varepsilon_i x_i}{a_i}=1,
 \qquad \varepsilon_i\in\{-1,1\}.
\end{align*}
Their feet are $(\varepsilon_1/a_1,\ldots,\varepsilon_d/a_d)/S$, so
\begin{align*}
 \D^2(B)=\prod_{i=1}^{d}\left[-\frac1{Sa_i},\frac1{Sa_i}\right].
\end{align*}
For~this new box the~sum of~the~inverse squared half-lengths is $S^2T$. Applying the~same formula again gives half-lengths $a_i/(ST)$.
\end{proof}

Thus planar polygons have a~universal return governed by~their number of~vertices, whereas spatial returns can depend on~the~geometry of~the~polytope and on~the~chosen point. Concentric spheres and centred boxes give two elementary ways in~which an~exact return survives.

\section*{Statements and Declarations}

\noindent\textbf{Funding.} The~author received no financial support for the~research, authorship, or publication of this article.

\noindent\textbf{Competing interests.}
The~author declares that he has no competing interests.

\noindent\textbf{Data availability.} No datasets were generated or analysed in this theoretical study.

\noindent\textbf{Use of generative artificial intelligence.} During the~preparation of this manuscript, the~author used OpenAI's ChatGPT as an auxiliary tool to assist with literature search, language editing and checking mathematical computations. All AI-assisted content was independently reviewed and verified by the~author, who takes full responsibility for the~manuscript.


\begin{thebibliography}{99}
\bibitem{Casey}
J.~Casey, \emph{A~Sequel to~the~First Six Books of~the~Elements of~Euclid}, 6th~ed., revised and enlarged, Hodges, Figgis \& Co., Dublin, 1892.
\bibitem{CoxeterGreitzer}
H.~S.~M.~Coxeter and S.~L.~Greitzer, \emph{Geometry Revisited}, New Mathematical Library, vol.~19, Mathematical Association of~America, Washington, DC, 1967.
\bibitem{deVilliers}
M.~de~Villiers, \emph{From nested Miquel triangles to~Miquel distances}, Math. Gaz. \textbf{86} (2002), no.~507, 390--395, \doi{10.2307/3621130}.
\bibitem{HumenbergerDeVilliers}
H.~Humenberger and M.~de~Villiers, \emph{Optimising Miquel circles centre triangles and third pedal triangles}, Math. Gaz. \textbf{110} (2026), no.~578, 287--304, \doi{10.1080/00255572.2025.2539599}.
\bibitem{Kostanyan}
A.~Kostanyan, \emph{Pedal tetrahedrons}, Mathematics in~School \textbf{6} (2026), no.~119, 92--99 (in~Armenian, with Russian and English summaries), \url{https://mathinschool.aspu.am/index.php/mathinschool/article/view/58}.
\bibitem{McDonaldGarciaReznik}
N.~McDonald, R.~Garcia and D.~Reznik, \emph{Exploring the~dynamics of~the~circumcenter map}, Math. Intelligencer \textbf{45} (2023), 232--241, \doi{10.1007/s00283-022-10223-9}.
\bibitem{Stewart}
B.~M.~Stewart, \emph{Cyclic properties of~Miquel polygons}, Amer. Math. Monthly \textbf{47} (1940), no.~7, 462--466, \doi{10.1080/00029890.1940.11991001}.
\bibitem{Thebault1914}
V.~Th\'ebault, \emph{Sur les polygones podaires d'un polygone donn\'e}, Nouv. Ann. Math. (4) \textbf{14} (1914), 433--446, \href{https://www.numdam.org/item/NAM_1914_4_14__433_0/}{Numdam}.
\bibitem{Thebault1941}
V.~Th\'ebault, \emph{Bibliography on~cyclic properties of~Miquel polygons}, Amer. Math. Monthly \textbf{48} (1941), 541, \doi{10.2307/2303391}.
\bibitem{Ziegler}
G.~M.~Ziegler, \emph{Lectures on~Polytopes}, Graduate Texts in~Mathematics, vol.~152, Springer, New York, 1995, \doi{10.1007/978-1-4613-8431-1}.
\end{thebibliography}
\end{document}